\documentclass[12pt,leqno,twoside]{amsart}
\usepackage{amssymb}
\usepackage{amsmath,amsfonts}
\usepackage{amscd}
\usepackage{latexsym}
\usepackage{verbatim}
\usepackage{epsfig}

\newtheorem{theorem}{Theorem}[section]

\theoremstyle{definition}

\theoremstyle{remark}

\theoremstyle{remark}

\theoremstyle{remark}

\theoremstyle{remark}

\theoremstyle{remark}

\theoremstyle{remark}

\renewcommand{\Bbb}{\mathbb}
\newcommand{\cal}{\mathcal}

\newcommand{\bC}{{\Bbb C}}

\def\growarrow#1{
  \setbox1=\hbox{ $\scriptstyle #1$\ }

\mathop{\smash{\hbox to \wd1{\rightarrowfill}}
          \vphantom\rightarrow}\limits^{#1}}

\begin{document}

\title[A supplement to a theorem of Merker and Porten: A short proof
of Hartogs' extension theorem for $(n-1)$-complete spaces]
{A supplement to a theorem of Merker and Porten: A short proof 
of Hartogs' extension theorem for $(n-1)$-complete complex spaces}
\thanks{Submitted on November $12^{th}$ 2008}

\author{Mihnea Col\c toiu}
\address{{\rm M. Col\c toiu}:
Institute of Mathematics of the Romanian Academy, P.O. Box 1-764, RO-014700, Bucure\c sti, Romania.}
\email{Mihnea.Coltoiu@imar.ro}





\maketitle

\section{Introduction}

The well-known Hartogs' extension theorem states that for every open subset
$D \subset \bC^n$, $n \geq 2$, every compact subset $K \subset D$ such that
$D \setminus K$ is connected the holomorphic functions on $D \setminus K$
extend to holomorphic functions on $D$.
For a simple and short $\bar \partial$ proof see \cite{E}. A long and very
involved proof of this result on 33 pp. ,using Morse theory, in the spirit
of Hartogs' original idea \cite{Ha} of  moving discs to get the extension,was
recently obtained by J. Merker and E. Porten \cite{M-P1}.

The Hartogs' theorem was generalized to $(n-1)$-complete manifolds (in the
sense of A. Andreotti and H. Grauert \cite{A-G}) by A. Andreotti and D. Hill
\cite{A-H} using cohomological results ($\bar \partial$ method).
In their forthcoming paper  \cite{M-P2} J. Merker  and E. Porten observed that in the singular case ``it is at present
advisable to look for methods avoiding $\bar \partial$ methods, because such
tools are not yet available'' and ``the essence of the present article is to
transfer such an approach to $(n-1)$-complete general complex spaces, where
the $\bar \partial$ techiques are still lacking, with some new difficulties due to singularities''.
J. Ruppenthal \cite{R} developped a $\bar \partial$ machinery for proving the
Hartogs' extension  theorem on Stein spaces with isolated singularities.

The main result of J. Merker and E. Porten \cite{M-P2}, which generalizes
Andreotti-Hill theorem \cite{A-H} for singular spaces, can be stated as follows:

\begin{theorem}

Let $X$ be a normal $(n-1)$-complete space ($n=dim X$), $D \subset \subset X$
a relatively compact open subset, $K \subset D$ a compact subset such that
$D \setminus K$ is connected. Then every holomorphic function on $D \setminus K$ can be extended to a holomorphic function on $D$.

\end{theorem}

In fact they proved this result even for the extension of meromorphic functions
(previously considered in the smooth case by V. Koziarz and F. Sarkis \cite{K-S}) but we
shall consider in this short note only the holomorphic extension.
The 20 pages proof of Merker and Porten \cite{M-P2} is also based on their previous
paper \cite{M-P1} on 33 pp., so putting together one gets about 50 pages which are very
technical.

We will give in this short note a 1 page proof for Theorem 1.1., using the
$\bar \partial$ method on the resolution of singularities, especially the Takegoshi
relative vanishing theorem \cite{T}, ( see also \cite{O}), which gives even
a more general statement valid on cohomologically $(n-1)$-complete spaces
( and without the assumption that $D$ is relatively compact).
Namely one has:

\begin{theorem}

Let $X$ be a $n$-dimensional normal cohomologically $(n-1)$-complete
complex space, $D \subset X$ an open subset, $K \subset D$ a compact
subset such that $D \setminus K$ is connected.Then every holomorphic function on $D \setminus K$ can be extended to a holomorphic function on $D$.

\end{theorem}
The ideas of the proof are essentially contained in the paper \cite{C-S}

\section{Proof of the result}

For the basic definitions of $q$-convex functions, $q$-complete complex space we
reffer to \cite{A-G}. We also recall that a complex space $X$ is called cohomologically
$q$-complete if one has the vanishing of the cohomology groups $H^i(X,\cal F)=0$ for every $i \geq q$ and every $\cal F \in Coh(X)$. By the main result
of \cite{A-G} a $q$-complete space is cohomologically $q$-complete ( a 
counter-example to the converse is still unknown).
For a complex manifold $X$ we denote by $K_X$ its canonical sheaf ( associated
to the canonical line bundle). Let $X$ be a complex (reduced) space and $\pi:
\tilde X \to X$ a resolution of singularities ( which exists by \cite{A-H-V},
\cite{B-M}). The following result , due to K. Takegoshi \cite{T} (see also T. Ohsawa \cite{O}), 
will be fundamental for our proof:

\begin{theorem}

Let $\pi:\tilde X \to X$ be a resolution of singularities of a complex space
$X$. Then one has the following vanishing for the higher direct images:
$R^i\pi_{*}K_{\tilde X}=0$ if $i \geq 1$

\end{theorem}

Let us also recall that by Grauert's coherence theorem \cite{G}
 $\pi_{*}K_{\tilde X}$ is a coherent sheaf on $X$. If moreover $X$ is assumed to be cohomologically $(n-1)$-complete it then 
follows that $H^i(X, \pi_{*}K_{\tilde X})=0$ if
$i \geq n-1$. By  Theorem 2.1.  the maps $H^i(X, \pi_{*}K_{\tilde X}) \to 
H^i(\tilde X, K_{\tilde X})$ are isomorphisms, so that one gets the vanishing
of the cohomology group
$H^i(\tilde X,K_{\tilde X})=0$ if $i \geq n-1$. By Serre duality \cite{S}
one gets the vanishing of the first cohomology group with compact supports
$H^1_c(\tilde X, \cal O_{\tilde X})=0$.
Moreover the arguments of L. Ehrenpreis \cite{E} ( see also \cite{Ho}) show
without any modification that the following holds: If $\tilde X$ is a
complex connected  non-compact manifold such that $H^1_c(\tilde X,
\cal O_{\tilde X})=0$ then for every open subset $\tilde D \subset \tilde X$
and for every compact subset $\tilde K \subset \tilde D$, such that
$\tilde D \setminus \tilde K$ is connected, it follows that every holomorphic
function on $\tilde D \setminus \tilde K$ can be extended to a holomorphic
function on $\tilde D$. Applying this result to $\tilde D=\pi^{-1}(D)$ and
$\tilde K=\pi^{-1}(K)$, where $\pi:\tilde X \to X$ is a resolution of
singularities for $X$, one gets immediately Theorem 1.2., since $\pi_{*}
\cal O_{\tilde X}=\cal O_X$ ( by the normality of $X$).



\begin{thebibliography}{MM}
\footnotesize{

\bibitem[A-G]{A-G}
A. Andreotti and H. Grauert: {\em Th\'eor\`eme de finitude pour la cohomologie
des espaces complexes}, Bull. Soc. Math. France 90 (1962), 193-259.

\bibitem[A-H]{A-H}
A. Andreotti and C.D. Hill: {\em E.E. Levi convexity and the Hans Lewy problem I and II}, Ann. Sc. Norm. Sup. Pisa 26 (1972), 325-363, 747-806.

\bibitem[A-H-V]{A-H-V}
J.M. Aroca, H. Hironaka, J.L. Vincente : {\em Desingularization theorems},
Mem. Math. Jorge Juan, No. 30, Madrid, 1977.

\bibitem[B-M]{B-M}
E. Bierstone and P. Milman: {\em Canonical desingularisation in characteristic
zero by blowing up the maximum strata of a local invariant}, Invent. Math.
128(1997), no. 2, 207-302.

\bibitem[C-S]{C-S}
M. Coltoiu and A. Silva : {\em Behnke-Stein theorem on complex spaces with singularities}. Nagoya J. Math.
137 (1995), 183-194.

\bibitem[E]{E}
L. Ehrenpreis: {\em A new proof and an extension of Hartogs' extension theorem}, Bull. Amer. Math. Soc. 67 (1961), 507-509.

\bibitem[G]{G}
H. Grauert: {\em Ein Theorem der analytischen Garbentheorie und die
Modulr\"aume komplexer Structuren}, Inst. Hautes \'Etudes Sci. Publ. Math.
no. 5 , 1960, 64 pp.

\bibitem[Ha]{Ha}
F. Hartogs : {\em Zur theorie der analytischen Funktionen mehrerer unabh\"angiger Ver\"anderlichen, unbesondere \"uber die Darstellung derselben durch Reichen,
welche nach Potenzen einer Ver\"anderlichen fortschreiten}, Math. Ann. 62 (1906), no. 1, 1-88.

\bibitem[Ho]{Ho}
L. H\"ormander : {\em An introduction to complex analysis in several
variables}, D. Van Nostrand Co., Inc., Princeton, N.J.-Toronto,1966, 208 pp.

\bibitem[K-S]{K-S}
V. Koziarz and F. Sarkis: {\em Probl\`eme du bord dand les vari\'et\'es
$q$-convexes et ph\'enom\`en de Hartogs-Bochner}, Math. Ann. 321 (2001),
no.3, 569-585.

\bibitem[M-P1]{M-P1}
J. Merker and E. Porten: {\em A Morse theoretical proof of the Hartogs
extension theorem}, J. Geom. Anal. 17 (2007), no.3, 513-546.

\bibitem[M-P2]{M-P2}
J. Marker and E. Porten: {\em The Hartogs' extension theorem on $(n-1)$-complete complex spaces}, preprint, arXiv:0704.3216 ( to appear in J. reine angew. Math).

\bibitem[O]{O}
T. Ohsawa : {\em A vanishing theorem for proper direct images}, Publ. RIMS
23(1987),no. 2, 243-250.

\bibitem[R]{R}
J. Ruppenthal : {\em A $\bar \partial$ theoretical proof of Hartogs' extension
theorem on Stein spaces with isolated singularities}, J. Geom. Anal. 18 (2008),no.4

\bibitem[S]{S}
J. P. Serre : {\em Un th\'eor\`eme de dualit\'e}, Comment. Math. Helv. 29(1955),9-26.

\bibitem[T]{T}
K. Takegoshi: {\em Relative vanishing theorems in analytic spaces}, Duke Math. J. 52 (1985), no. 1, 273-279.
 

}

\end{thebibliography}
\end{document}